\documentclass{article}
\usepackage{graphicx} 
\usepackage{amsfonts}
\usepackage{amsmath}
\usepackage{amssymb}
\usepackage{amsthm}
\newcommand{\intl}[2]{\operatornamewithlimits{\int}\limits_{#1}^{#2}}
\newcommand{\thetaf}[2]{\theta(\alpha+\intl{#1}{#2}\overrightarrow{\omega})}
\newcommand{\suml}[2]{\operatornamewithlimits{\sum}\limits_{#1}^{#2}}
\newcommand{\psl}[0]{PSL(2,\mathbb C)}

\title{Fourier decay, Green's functions and Schottky groups}
\author{Ilyas Bayramov}
\date{August 2023}

\begin{document}

\maketitle

\begin{abstract}
    In this note, I would like to discuss an approach to the construction of Green's function of laplacians on algebraic surfaces, indicated by Manin in \cite{manin}, towards the computation
    of the Green's function on surfaces using Schottky uniformization. We shall see that the exact geometric interpretation of the formula mentioned there is obscure, and
    try to remedy the situation by investigating convergence of
    deformations of that formula.
\end{abstract}

\section{Acknowledgments}

I thank Peter Kosenko and Dennis Sullivan for many useful suggestions.

\section{Motivation for Manin's formula}

There is a certain generalization due to Arakelov of the usual intersection pairing between the divisors that is useful, for example, in the proof of Mordell's theorem. Here by divisor we mean a Weil divisor, i.e., a complex linear combination of points on the Riemann surface $S$.

To each divisor $a$ of degree $0$ on a (smooth) algebraic curve $S$, by Riemann-Roch \cite{gross}, corresponds a meromorphic differential form $\omega_a$ such that $$a=\suml{x\in S}{}res_x(\omega_a)\cdot x.$$ 

This differential form is unique up to an addition of a holomorphic differential form to it. Since the real parts of periods of a holomorphic differential determine it uniquely (by, say, Riemann period relations) and we have the freedom to choose any holomorphic differential, we can insist that $\omega_a$ has all of its periods purely imaginary. 

Then $\text{d}g_a:=\omega_a-\bar\omega_a$ is exact since all of its periods are $0.$ $g_a$ is then called the Green-Arakelov function.

\textbf{Definition}\cite{gross}: Let $b=\operatornamewithlimits{\sum}\limits_ym_y(y),$ a divisor of degree $0$; let $a$ be another such divisor, $|a|\cap|b|=\emptyset$, $(y)\in S,$ $S$ an algebraic curve. Then the \textbf{height pairing} between $a$ and $b$ is $g(a,b)=\operatornamewithlimits{\sum}\limits_ym_yg_a(y).$

The goal is to compute $g_a.$ To do that, we need to take a detour into low-dimensional topology.

\section{Schottky uniformization: existence}

Schottky uniformization of a Riemann surface $S$ is a representation $\rho: \pi_1(S) \to PSL(2, \mathbb C),$ satiasfying the following conditions. Its image, called a Schottky group, is a free group $\Gamma$ that is purely hyperbolic, i.e., the absolute values of traces of its elements are strictly bounded by $2$ from below. $\psl$ acts on $\mathbb H^3$ by isometries, and on $S^2$ by linear fractional transformations. $S^2$ can be identified with the boundary of $\mathbb H^3$ at infinity (aka "the conformal infinity"). 

Now, any infinite discrete group acting on a compact space, in this case $\mathbb H^3 \cup S^2,$ must have accumulation points. The set of such accumulation points must then be contained in $S^2$ (because $\psl$ acts on $\mathbb H^3$ by isometries). For a given $\Gamma$ we denote it by $\Lambda(\Gamma)$ and call it the limit set. Its complement in $\S^2$ is denoted by $\Omega(\Gamma)$ and is called the domain of discontinuity. By Sullivan's dictionary, these correspond to Julia and Fatou sets, respectively.

"Uniformization" then comes from the application of Ahlfors' finiteness theorem.

It states that for any finitely generated discrete subgroup  of $\psl$, $$\Omega(\Gamma)/\Gamma$$ is a Riemann surface of finite area. Hence, in this particular case it is, as well, and to each Riemann surface we can associate a corresponding Schottky uniformization by Koebe's retrosection theorem.

\section{Schottky uniformization: uniqueness}

This uniformization is not unique; to see where exactly the non-uniqueness arises from, we must consider the action of a Schottky group $\Gamma$ on all of $\mathbb H^3\cup \Omega(\Gamma).$ Then the quotient $K(\Gamma)$ of such an action admits a hyperbolic metric in its interior, and it is actually evident that $K(\Gamma)$ is a handlebody. Then $S:=\Omega(\Gamma)/\Gamma$ determines a point in the Teichmüller space up to the action of the elements of the mapping class group that have representatives extending to homeomorphisms of $K(\Gamma)$ that are isotopic to the identity. This subgroup of $\text{Mod}(S)$ is denoted by $\text{Mod}_0(S)$ and is known as the handlebody group of a fixed genus $g.$

\section{Compactification $\bar{\mathfrak T}$ of the Teichmuller space}

Thurston defined a compactification of the Teichmuller space as follows: consider the set $\mathcal{S}$ of isotopy classes of simple closed curves on a surface. Then for two such classes $a,b\in \mathcal S$ the geometric intersection number $i(a,b)$ is given by the smallest number of intersections between the representatives of $a$ and $b$, respectively. Denote by $\mathbb R_+$ the set of nonegative reals with the usual induced topology. Then there exists a map $$i_{*}: \mathcal S\to \mathbb R_+^{\mathcal S}\backslash 0,$$ $$i_{*}(a)(b)=i(a,b).$$ Now, there also exists a map $\pi: \mathbb R_+^{\mathcal{S}}\backslash 0\to P(\mathbb R_+{\mathcal{S}}).$ Here arises the question of topology on these spaces: we take $\mathbb R_+^{\mathcal S}$ with the weak (product) topology, and $P(\mathbb R_+^{\mathcal S})$ with the quotient topology. Then the closure $\bar{\mathcal S}$ of $\pi\circ i_*(\mathcal S)$ in quotient topology is homeomorphic to a sphere $S^{6g-6-2b},$ where $b$ is the number of punctures. 

On the other hand, there also exists an embedding of Teichmuller space into $P(\mathbb R^{\mathcal S}).$ It is given by the length function $$l: \mathfrak T \times \mathcal S\to \mathbb R_+^{\mathcal S},$$ $$l(\theta, \alpha)=\inf\limits_{\gamma\in \alpha}\theta(\gamma),$$ where $\gamma$ is a simple closed curve in the isotopy class of $\alpha.$ Since this is well-defined, because both $\theta$ and $\alpha$ are defined up to isotopy, this gives an embedding $l_*:\mathfrak T\to \mathbb R_+^{\mathcal S}\backslash 0$ which is a homeomorphism onto the image. 

Now it is obvious that the geometric intersection number of $a\in \mathcal S$ with itself is $0,$ and by \cite{bonahon}, the geometric intersection number can be extended to $\mathfrak T$ and $\bar{\mathcal S}$ in such a way that the self-intersection number is respectively a non-zero constant on all of $l_*(\mathfrak T)$ and $0$ on $\bar{\mathcal S}.$ Since we consider the weak topology on $\mathbb R_+^{\mathcal S}\backslash 0,$ we see that if a sequence $\theta_n$ of elements of $l_*(\mathfrak T)$ leaves any compact set of $\mathbb R_+^{\mathcal S}\backslash 0$ in such topology, then we must have that for at least one of the elements $\alpha$ of $\mathcal S$ $l(\theta_n, \alpha)\to \infty.$ Then if a sequence $\pi(\theta_n)$ converges to some point $b$ in $P(\mathbb R^{\mathcal S}),$ a sequence of their preimages in $\mathbb R_+^{\mathcal{S}}$ must converges to a function $a\in \pi^{-1}(b);$ that means that there exists numbers $\lambda_n>0$ such that $\lambda_n\theta_n$ converges to $a;$ as per the above, $\lambda_n\to 0,$ as at least one of the coordinates $(\theta_n(s))_{s\in \mathcal S}$ goes to $\infty.$ Then by the bilinearity (!) of the intersection number $i(\lambda_n\theta_n,\lambda_n\theta_n)=|\lambda_n|^2i(\theta_n,\theta_n)=C|\lambda_n|^2\to 0=i(a,a).$ This implies that $a\in \bar {\mathcal S},$ by the converse of the above statement about $\bar {\mathcal S},$ which is less obvious. Notice that this argument is very similar to a measure concentration phenomenon.
Thus, $\bar{\mathfrak T} = \pi\circ l_*(\mathfrak T)\cup {\mathcal S}$ is a compactification of $\mathfrak T.$

\section{The action of $\text{Mod}_0(S)$ on $\bar{\mathfrak T}$}

Profoundly, the situation in this case is reminiscent of the Kleinian group action on $\mathbb H^3:$ for the subgroup $\text{Mod}(H)$ of $\text{Mod}(S)$ that consists of elements that have representatives extending to diffeomorphisms of the handlebody $H,\partial H=S$, there is a domain of discontinuity on $\bar{\mathfrak T},$ and a limit set that is contained in $\bar{\mathcal S}.$ Now, following Masur \cite{masur}, lemma 1.1, we have that if a curve on $S$ bounds a disk in the corresponding handlebody, and $\alpha_1, \dots, \alpha_g$ are simple closed curves that bound disks in $H$ and cut $S$ into a $2g$ holed sphere, then $\beta \in \mathcal{S}$ bounds a disk iff $i(\beta,\alpha_i) = 0$ or we have the following situation. Suppose $i(\beta, \alpha_i)\neq 0$ for some $i$. Choose $(\beta_l)_{l\in \{1,\dots, g\}}\in \mathcal S$ with $i(\beta_l, \alpha_j)=\delta_{lj}$ and $\beta_l-$free generators of $\pi_1(H)=\Gamma.$ As a word in $\beta_l$ and $\beta_l^{-1},\beta$ is trivial. Up to a cyclic permutation, there must be either $\beta_j\beta_j^{-1}$ or $\beta_j^{-1}\beta_j$ somewhere in the word. But then $\beta$ crosses $\alpha_j$ and then crosses it in the opposite direction.

Now, consider $\mathcal S_c\subset \mathcal S,$ the set of simple closed curves bounding disks in the handlebody. Then by Theorem 2.2 in Masur $\bar{\mathcal S_c}$ is the limit of the actions of both $\text{Mod}(H)$ and $\text{Mod}_0(S).$

Finally, there is the canonical measure $m$ on $\bar{\mathcal S},$ see Masur's theorem $B,$ such that $m(N:=\{a\in \bar{\mathcal S}-\bar{\mathcal S_c}|i(a,b)=0 \text{ for some }b\in \bar{\mathcal S_c}\}) = 0.$ Taking $\mathcal O := (\bar{\mathcal S}-\bar{\mathcal S_c}) - N,$ we have that $\mathcal O$ is the domain of discontinuity of the action of both $\text{Mod}(H)$ and $\text{Mod}_0(S)$ on $\mathcal S$ by theorem 2.1 in Masur. Moreover, in genus $2,$ according to Masur, the limit set $\bar{\mathcal S_c}$ has measure $0.$

\section{Deforming the Schottky group}

We consider here $AH(H)$ to be the set of all discrete faithful representations of $\Gamma$ into $\psl$ up to conjugation. Now, according to Li-Ohshika-Lecuire \cite{ohshika}, we have that if $\rho_n\in AH(H)$ corresponds to a sequence of elements of $\mathfrak T$ that converges to an element of Masur domain, then $\rho_n$ converges to an element of $AH(H).$ Finally, we have, according to Anderson-Canary-Culler-Shalen \cite{accs}, that the set of Kleinian groups with the limit set the whole sphere is dense in the boundary of the domain of convex-cocompact representations, i.e., in this case, the quasi-conformal deformations of a Schottky group. Thus, if a sequence of elements $a_n$ of $\mathfrak T$ leaves any compact subset of $\mathfrak T$ in the quotient of weak topology described above, then its limit $a\in \mathfrak T$ almost always corresponds to a group $\Gamma_a$ with limit set $\Lambda(\Gamma_a)=S^2.$ This is why we can \textbf{almost} quasiconformally deform a Schottky group to another one that has $S^2$ as its limit set.

\section{Manin's formula}

In what follows, $\langle a,b,c,d\rangle$ denotes the cross-ratio of $a,b,c,d.$

Using Schottky uniformization, Manin \cite{manin} defined $\omega_a$ for any divisor of degree $0$, as in section $1$; for $a=(x)-(y)$ it is $$\nu_{(x)-(y)} - \operatornamewithlimits{\sum}\limits_lX(x,y)\omega_l,$$ where $$\nu_{(x)-(y)}=\text{dlog}W_{(x)-(y),z_0}(z)$$, 
$$X_l(a,b)= -\sqrt{-1}\frac{\operatornamewithlimits{\sum}\limits_{h\in S(g_k)}\text{log}|\langle a,b,z^+(h),z^-(h)\rangle|}{\operatornamewithlimits{\sum}\limits_{h\in C(g_k|g_l)}\text{log}| \langle z^+(g_k)z^-(g_k),hz^+(g_l),hz^-(g_l)\rangle |},$$ $$ \omega_l = -\sqrt{-1}\frac{1}{2\pi}\text{dlog}W_{(g_kz_1)-(z_1),z_0}(z).$$ Here, $$W_{d,q}(r)=\operatornamewithlimits{\prod}\limits_{h\in\Gamma}\frac{w_d(hr)}{w_d(hq)},$$ for $d$ a divisor in  $\mathbb{C}P^1,$ $w_d$ is the meromorphic function with this divisor. Also, $S(g_k)$ is the conjugacy class of $g,$ and $C(g_k|g_l)-$double coset representatives.

\section{Simplification of Manin's formula}

It is known as a consequence of Fay's trisecant theorem \cite{poor} that $$\text{exp}(\intl{a_1}{a_2}\omega_{(x)-(y)})=\frac{\thetaf{a_1}{x} \thetaf{a_2}{y}}{\thetaf{a_2}{x}\thetaf{a_1}{y}},$$ $$Re\intl{a_1}{a_2}\omega_{(x)-(y)} = \log(\thetaf{a_1}{x}) +\log(\thetaf{a_2}{y})$$ $$-\log(\thetaf{a_2}{x})-\log(\thetaf{a_1}{y});$$
then by Roelcke and Fay \cite{kokotov} we have $$\log(\thetaf{x}{y})=\int{dS(p)}\int{dS(q)}Re\intl{x}{p}\omega_{(y)-(q)},$$ where $$Re\intl{x}{p}\omega_{(y)-(q)}=\operatornamewithlimits{\sum}\limits_{h\in \Gamma}\log|\langle y, q, hx, hp\rangle|;$$ so we have as the Green's function $$\intl{\Omega(\Gamma)}{}\intl{\Omega(\Gamma)}{}Re\intl{x}{p}\omega_{(y)-(q)}dS(p)dS(q)$$
where $dS$ is the $\Gamma-$invariant measure on $\Omega(\Gamma).$

\section{Fixed-point theorems}

We would like to rephrase the existence of the solution to Laplace's equation on $\Lambda(\Gamma)$ with potential (which is equivalent to being able to find Green's function) through a fixed point problem. We follow Evans \cite{evans} in this. Consider the following setup: assume that the smooth functions $\{w_k\}_{k\in \mathbb N}$ form an orthonormal basis of $W^{1,2}_0$, where $W^{1,2}_0$ is the closure of the space of compactly supported smooth functions in $W^{1,2},$ the Sobolev space of square-integrable functions with square-integrable weak derivative. For example, we can take $w_k$ to be harmonics so that the corresponding coefficients of $w_k$ in a series converging to the function $f\in W^{1,2}_0$ are Foruier coefficients. Consider the system of equations $$\intl{U}{}\mathbf{a}(Du_m)\cdot Dw_kdx = \intl{U}{} fw_kdx,$$ where $u_m=\suml{k=1}{m} d^k_mw_k,$ but for each $u_m$ we only consider those equations up to $k=m.$ Here $\mathbf{a}:\mathbb R^n\to \mathbb R^n, \mathbf{a}=\{a^1,\dots, a^n\}$ is a vector field. Also, we must impose the following on $\mathbf a:$ $$(\mathbf{a}(p)-\mathbf{a}(q))(p-q)\ge 0,$$ $$\mathbf{a}\le C(1+|p|),$$ $$\mathbf a(p)\cdot p\ge \alpha |p|^2 - \beta$$ for some $C,\alpha>0$ and $\beta\ge 0.$ Then we can use the fixed point theorem to show that, if we consider $$v^k(d) := \intl{U}{}\mathbf{a}(\suml{j=1}{m} d_jDw_j) \cdot Dw_k - fw_kdx (k=1,\dots,m),$$ over all  $d=(d_1,\dots, d_m)\in \mathbb R^m,$ then $v^k(d)$ has a zero. 

This equation is equivalent to Laplace's equation $-\Delta u + cu = 0,$ where $c =\frac{2a\Delta a-(\nabla a)^2}{4a^2}$ (see \cite{div}, for example).

We can interpret the space of $d'$s in this case as a space of Fourier coefficients. Naud et al. \cite{lnp} have established (non-uniform) boundedness of Fourier coefficients of functions on the limit set of Schottky groups. Thus, we have that the Green's function exists on $\Lambda(\Gamma)$ for any $\Gamma$ free convex-cocompact.

Now, by \cite{mjseries}, we have that the limit sets of algebraically convergent sequence of Schottky groups, defined in section \textbf{6}, converge to the limit set of the of the group $\Gamma_a$ in Gromov-Hausdorff sense. On the other hand, following Quint \cite{quint}, we have that the Patterson-Sullivan measure is in the Lebesgue class of measures on $S^2,$ hence it is doubling with respect to the Fubini-Study metric on $S^2,$ which we identify with $\mathbb CP^1$. By the results of \cite{ding} adapted to the case of sets, rather than manifolds, we have that if we have a sequence of metric spaces equipped with measures that are doubling with respect to those metrics, then the Green's functions on those sets converge to a Green's function on the Gromov-Hausdorff limit of those sets.

Finally, considering both $$\operatornamewithlimits{\sum}\limits_{h\in \Gamma_a}\log|\langle y, q, hx, hp\rangle|$$ and $$\operatornamewithlimits{\sum}\limits_{h\in \Gamma_{a_n}}\log|\langle y, q, hx, hp\rangle|,$$ $\Gamma_a$ being the algebraic limit of $\Gamma_{a_n}$ as above, we have $$\operatornamewithlimits{\sum}\limits_{h\in \Gamma_{a}}\log|\langle y, q, hx, hp\rangle|=$$$$\operatornamewithlimits{\sum}\limits_{h\in \Gamma_{a}}\chi_{\Omega(\Gamma_{a_n})}(\langle y, q, hx, hp\rangle)\log|\langle y, q, hx, hp\rangle|+$$ $$+\operatornamewithlimits{\sum}\limits_{h\in \Gamma_{a}}\chi_{\Lambda(\Gamma_{a_n})}(\langle y, q, hx, hp\rangle)\log|\langle y, q, hx, hp\rangle|,$$ where $\chi'\text{s}$ are indicator functions. Since the action of $\Gamma_a$ is topologically dense on $S^2,$ that is, the closure of any orbit is dense, we have that, in particular, any orbit of $\Gamma_{a_n}$ that lies in $\Omega(\Gamma_{a_n})$ can be approximated by the orbit of the same point under $\Gamma_a.$ Since the action of $\Gamma_{a_n}$ on $\Lambda(\Gamma_{a_n})$ is also topologically dense, we can estimate $$\operatornamewithlimits{\sum}\limits_{h\in \Gamma_{a}}\chi_{\Lambda(\Gamma_{a_n})}(\langle y, q, hx, hp\rangle)\log|\langle y, q, hx, hp\rangle|$$ by $$\operatornamewithlimits{\sum}\limits_{h\in \Gamma_{a_n}}\chi_{\Lambda(\Gamma_{a_n})}(\langle y, q, hx, hp\rangle)\log|\langle y, q, hx, hp\rangle|,$$ which converges by \cite{lnp}.

Thus, we have the following:

\textbf{Theorem:} For almost any Schottky representation $\Gamma_0$ of a fixed free group, for almost any point at the boundary of the Schottky space, there exists an unbounded sequence of representations $\Gamma_n$ converging to it and quasiconformal maps on $S^2_\infty,$ conjugating the action of $\Gamma_0$ to $\Gamma_n,$ such that for any $n$, the expression $\operatornamewithlimits{\sum}\limits_{h\in \Gamma_{a_n}}\log|\langle y, q, hx, hp\rangle|$, where $\Gamma$ is the limit of $\Gamma_n$ in $AH(S),$ for any quadruple of points $y,q,x,p$, converges.

\textbf{Proof:} If $\langle y, q, x, p\rangle\in \Lambda(\Gamma_{a_n}),$ this follows by \cite{lnp} and invariance of $\Lambda(\Gamma_{a_n})$. If  $\rangle y, q, x, p\langle\in \Omega(\Gamma_{a_n}),$ then we can estimate it from above by $$\operatornamewithlimits{\sum}\limits_{h\in \Gamma_{a}}\chi_{\Omega(\Gamma_{a_n})}(\langle y, q, hx, hp\rangle)\log|\langle y, q, hx, hp\rangle|,$$ which converges by the above. This concludes the proof.

\textbf{Remark:} The result above generalizes \cite{pollicott}.

\bibliographystyle{abbrv}
\bibliography{refs}

\end{document}